\documentclass{article}
\usepackage{amsmath}
\usepackage{amsfonts}
\usepackage{amssymb}
\usepackage{pstricks,pst-node,amssymb,amsmath,graphics,latexsym,tabularx,shapepar}

\usepackage{graphicx}
\usepackage{float}
\newtheorem{theorem}{Theorem}

\newcommand{\bpartial}{\mathop{\partial\kern -4pt\raisebox{.8pt}{$|$}}}
\newcommand{\bra}{\mathopen{[\kern-1.6pt[}}
\newcommand{\ket}{\mathclose{]\kern-1.5pt]}}
\newcommand{\bbra}{\mathopen{[\kern-2.2pt[\kern-2.3pt[}}
\newcommand{\bket}{\mathclose{]\kern-2.1pt]\kern-2.3pt]}}

\makeindex

\begin{document}


\title{Similarity Solutions for the Flux-limited Keller--Segel System with Time-Varying Chemical Decay
Rate}
\author{
Ahmed Abbas Jaber Al-Furaiji\thanks{Department of Mathematical Sciences, Azarbaijan Shahid Madani University, Tabriz, Iran. Email: 
ahmedgon1234567@gmail.com}
\and
Ghorbanali Haghighatdoost\thanks{Corresponder author, Department of Mathematical Sciences, Azarbaijan Shahid Madani University, Tabriz, Iran. Email: gorbanali@azaruniv.ac.ir} 
\and
Mustafa Bazghandi\thanks{Department of Mathematical Sciences, Azarbaijan Shahid Madani University, Tabriz, Iran. Email: mostafabazghandi2001@gmail.com}
}
\maketitle

\begin{abstract}
We investigate a one-dimensional flux-limited Keller--Segel system (FLKS) in which the chemical decay rate is allowed to vary explicitly in time, a feature motivated by enzymatic regulation and environmental variability in chemotactic signalling. Treating the decay rate as an arbitrary function, we perform a systematic Lie symmetry analysis and use equivalence transformations to classify the admitted symmetry extensions. The ordinary kernel algebra of the class is generated by spatial translations. In addition, each fixed decay function gives rise to a decay-dependent gauge symmetry in the chemoattractant variable, so that a generic individual system admits a two-dimensional Abelian Lie algebra. Group classification via the classifying equation shows this kernel extends only when decay function is constant. We construct an optimal system of subalgebras and derive the corresponding similarity reductions. Finally, we find some explicit solutions for our FLKS model. Our results provide a rigorous mathematical foundation for understanding which temporal decay patterns admit similarity reductions, thereby enabling analytical progress on flux-limited chemotaxis models with realistic time-varying degradation mechanisms.
\end{abstract}

\noindent{\bf AMS}: 35Q92, 70G65, 76M60.

\section{Introduction}\label{sec:introduction}
The classical Keller--Segel (KS) model, while foundational in the mathematical study of chemotaxis, exhibits a critical mathematical deficiency: finite-time blow-up of cell densities under supercritical initial mass conditions. Such singularities are biologically implausible, as real cell populations remain bounded by volume exclusion, nutrient depletion, and other regulatory mechanisms. Preventing or regularizing these blow-up phenomena has thus become a central objective in chemotaxis modeling. To this end, various extensions have been proposed--including logistic growth terms \cite{Xiang2018}, cross-diffusion mechanisms \cite{Carrillo2012,Hittmeir2011}, tempered fractional derivatives \cite{Haghighatdoost2025tempered}, and flux-limited formulations--each designed to restore global existence of solutions and produce physically meaningful long-time dynamics. Beyond the blow-up problem, the classical model suffers from an equally fundamental limitation in its representation of chemotactic transport. The classical model assumes velocity is linear ($\mathbf{u} = \chi \nabla v$), permitting unrestricted velocity growth as the magnitude of the gradient increases, failing to account for the inherent biomechanical constraints that govern cellular motility. 

\par 
Recent advances have demonstrated that flux limitation corrects several core shortcomings of classical chemotaxis models: it enforces finite propagation speed in contrast to the infinite-speed transport of the classical Keller--Segel system, and yields aggregation dynamics consistent with empirical chemotaxis observations. Current research provides both rigorous mathematical analysis and computational evidence supporting the FLKS framework as a practical and robust model for chemotactic transport \cite{Chertock2012,Erban2004,Jaiswal2023,Kohatsu2025,Kohatsu2025SelfSimilar,Painter2019}.

\par 
In many chemotactic systems, the dominant mechanism controlling signal degradation is enzymatic. For example, in \textit{Dictyostelium discoideum}, multiple phosphodiesterases regulate the hydrolysis of extracellular cAMP, and their expression levels vary substantially throughout the developmental cycle, causing strong temporal modulation of the effective decay rate $\kappa(t)$ \cite{Bader2007, Golding2005, Lenstra2016}. Similar forms of enzyme-mediated regulation of chemoattractant breakdown occur in bacterial chemotaxis, where ligand-specific periplasmic enzymes adjust signal lifetimes in response to environmental conditions \cite{Alber2019,Erban2004}. These biochemical processes typically operate on developmental or stress-induced timescales, and therefore naturally induce time-dependent degradation dynamics. Consequently, enzymatic regulation provides a mechanistic basis for modeling $\kappa$ as a non-constant function of time, highlighting the biological relevance of the generalized framework adopted in this study.

\par 
To accommodate these critical biological requirements, we analyze a generalized system where the chemoattractant degradation $\kappa$ is explicitly a function of time $t$. We restrict our analysis to the one-dimensional setting ($x \in \mathbb{R},\; t\geq 0$), a choice motivated by several complementary considerations. First, one-dimensional models capture universal scaling behavior and intermediate asymptotics that persist across dimensions, while isolating the essential competition between diffusion, chemotactic flux, and time-dependent decay without the geometric complexities inherent in radial symmetry or multi-directional gradients. Second, this restriction ensures analytical tractability: the reduced ordinary differential equations obtained through symmetry methods admit exact or qualitative analysis that is generally unattainable in higher dimensions, thereby providing rigorous benchmark solutions for validating numerical schemes in two and three dimensions. Finally, many biologically relevant chemotaxis scenarios occur naturally in quasi-one-dimensional geometries, including bacterial chemotaxis in capillary tubes and microfluidic channels \cite{Kalinin2009}, aggregation streams in Dictyostelium discoideum \cite{Palsson1996}, and cell migration along thin tissue strips or wound edges during tissue repair and morphogenesis \cite{Maini2004}.  

We consider the one-dimensional flux-limited Keller--Segel system
\begin{equation}
\label{eq:model-main}
\begin{cases}
u_t = D u_{xx} - \left( u F(v_x) \right)_x + r u\left(1-\dfrac{u}{K}\right), \\[2mm]
\tau v_t = v_{xx} - \kappa(t)v + u,
\end{cases}
\end{equation}
where
\begin{equation}
\label{eq:model-flux}
F(p) = \frac{p}{\sqrt{1+p^2}},
\end{equation}
is the algebraic limiter captures generic saturation of chemotactic velocity \cite{Perthame2019,Zeng2025},
The dependent variables $u=u(x,t)$ and $v=v(x,t)$ represent, respectively, the cell density and the concentration of the chemoattractant.  $D, \tau,r,K >0$ are constants. $\kappa(t)$ is an arbitrary function.

\par 
Lie symmetry analysis has become an essential tool in the study of nonlinear partial differential equations, as it enables a systematic identification of the intrinsic geometric structures underlying differential equations. These structures can then be exploited to analyse differential invariants, construct exact or similarity solutions, and reveal the associated symmetry properties of the system (e.g. \cite{Bazghandi2026,Haghighatdoost2025}). Complementing this framework, equivalence transformations play a crucial role when investigating families of equations containing arbitrary functions or parameters, since they allow one to map the class of equations to canonical representatives, eliminate redundant cases, and clarify which structural features are essential for determining the behaviour of the model (e.g. \cite{Gazizov2019}).

\par 
Our objective is to perform a systematic classification of Lie symmetries with respect to arbitrary decay functions $\kappa(t)$, employing equivalence transformation theory to identify all admissible symmetry algebras. The foundational theory of group classification for nonlinear equations was proposed by Ovsiannikov \cite{Ovsiannikov1959}, and has since been extended through both theoretical advances and widespread application to nonlinear evolution equations \cite{Ivanova2006,Ivanova2010,Opanasenko2017,Popovych2004}. This approach enables us to determine which temporal decay patterns admit similarity reductions, thereby providing analytical progress on flux-limited chemotaxis models with realistic time-varying degradation mechanisms. To the best of our knowledge, no previous Lie symmetry analysis of flux-limited Keller--Segel (FLKS) systems with time-varying chemical decay exists in the literature.

\par 
The remainder of this paper is organized as follows. In Section \ref{sec:preliminaries}, we introduce the equivalence transformation framework necessary for a systematic group classification of the class of systems parametrized by the arbitrary decay function $\kappa(t)$. Section \ref{sec:lie-analysis} derives the full set of Lie symmetry determining equations, identifying both the universal symmetries and the crucial classifying equation that governs possible extensions. The complete group classification is presented in Section \ref{sec:classification}, where we establish the kernel algebra for arbitrary decay rates and prove that the only symmetry extension occurs when $\kappa(t)$ is constant. Section \ref{sec:optimal-systems} constructs the optimal system of one-dimensional subalgebras, providing a non-redundant catalog of all admissible symmetry reductions. These reductions are then applied in Section \ref{sec:reductions} to derive invariant solutions, including novel spatially structured and traveling-wave solutions for the flux-limited system. We conclude in Section \ref{sec:conclusion} with a summary of our findings and potential directions for future research.


\section{Preliminaries}
\label{sec:preliminaries}

The classical Lie symmetry method is designed to determine the point symmetries of a single, fully specified differential equation in which every coefficient, parameter, and constitutive function is fixed in advance. However, system \eqref{eq:model-main} is not a single equation but rather a family of systems parametrized by the arbitrary decay function $\kappa(t)$. For such a class of equations, the standard symmetry analysis must be elevated to a group classification problem. The appropriate theoretical framework for this task is provided by the theory of equivalence transformations, which identifies admissible changes of variables that preserve the differential structure of the entire family while allowing the arbitrary element to transform.

Let $\mathcal{K}$ denote the class of flux-limited Keller--Segel systems \eqref{eq:model-main}. Each particular choice of $\kappa(t)$ determines one member of the class $\mathcal{K}$.

\subsection{Lie symmetries of a fixed member of the class}

For a prescribed decay function $\kappa(t)$, a Lie point symmetry is a local transformation of the variables
\[
(x,t,u,v)
\]
that maps solutions of the corresponding system \eqref{eq:model-main} to solutions of the same system. In this problem, $\kappa(t)$ is held fixed and is not transformed. The resulting maximal Lie invariance algebra will be denoted by $\mathcal A^{\max}(\kappa)$.

The intersection of these algebras over all admissible decay functions,
\[
\mathcal A^{\ker}=\bigcap_{\kappa}\mathcal A^{\max}(\kappa),
\]
is called the \emph{kernel algebra} of the class. It consists of the symmetries admitted for every choice of $\kappa(t)$. The remaining determining equations yield a classifying condition involving $\kappa$ and its derivatives. Solving this condition identifies the special forms of $\kappa(t)$ for which $\mathcal A^{\max}(\kappa)$ is larger than the kernel algebra.

\subsection{Equivalence transformations of the class}
Group classification of the class $\mathcal K$ requires more than solving the determining equations for each fixed $\kappa(t)$ separately: since a Lie point transformation between two systems in $\mathcal K$ that also reparametrizes $\kappa$ identifies their maximal invariance algebras $\mathcal A^{\max}(\kappa)$ and $\mathcal A^{\max}(\widetilde\kappa)$ as isomorphic (conjugate via that transformation), such pairs $\kappa,\widetilde\kappa$ do not need to be classified independently. This equivalence relation on the arbitrary element is generated by the \emph{equivalence group} of $\mathcal K$: the set of admissible point transformations of $(x,t,u,v)$, together with a corresponding transformation of $\kappa$, that map every member of \eqref{eq:model-main} to another member of the same class while preserving the fixed constants $D,\tau,r,K$ and the fixed flux function $F$. Determining this group is therefore the necessary first step of the classification: it fixes the equivalence classes of $\kappa(t)$ within which the subsequent Lie symmetry analysis (Section \ref{sec:lie-analysis}) and classification (Section \ref{sec:classification}) are to be carried out, and it supplies the canonical representatives and normalizations used throughout.

We seek the equivalence group in the standard projectable form
\[
\widetilde x=X(x),\qquad \widetilde t=T(t),\qquad
\widetilde u=U(u),\qquad \widetilde v=V(v),\qquad \widetilde\kappa=\widetilde\kappa(\widetilde t),
\]
with nonzero Jacobian, imposing that \eqref{eq:model-main} retain its functional form for some admissible $\widetilde\kappa$. Because the diffusion terms are fixed and linear, $X,T,U,V$ are necessarily affine; the fixed logistic term with fixed $r,K$ forces $U(u)=u$; the fixed, non-homogeneous, odd flux limiter $F$ then forces $X_x=\varepsilon\in\{-1,1\}$ (no admissible rescaling of $x$ or $v$); and the fixed diffusion and source coefficients together force $T_t=1$ and $V(v)=v$. The equivalence group of $\mathcal K$ is therefore
\begin{equation}
\label{eq:equivalence-group}
\widetilde x=\varepsilon x+\delta,\qquad
\widetilde t=t+\beta,\qquad
\widetilde u=u,\qquad
\widetilde v=v,
\qquad \varepsilon\in\{-1,1\},
\end{equation}
acting on the arbitrary element by
\begin{equation}
\label{eq:kappa-equivalence}
\widetilde\kappa(\widetilde t)=\kappa(\widetilde t-\beta).
\end{equation}
Thus $\mathcal K$ admits only space-time translations and spatial reflection ($\varepsilon=-1$, admissible because $F$ is odd) as equivalence transformations, with no scaling freedom in $x,t,u,v$ and no rescaling of $\kappa$. In particular, \eqref{eq:kappa-equivalence} only shifts the argument of $\kappa$, so a nonzero constant $\kappa\equiv\kappa_0$ can never be mapped to a different constant; distinct constants remain inequivalent representatives within $\mathcal K$, and this normalization is used directly in the classification of Section \ref{sec:classification}.


\section{Lie symmetry determining equations}
\label{sec:lie-analysis}
We now determine the Lie point symmetries of a fixed system \eqref{eq:model-main}. Thus, the function $\kappa(t)$ is no longer transformed.

\subsection{Infinitesimal generator}

Let
\[
\mathbf{X} = \xi(x,t,u,v)\partial_t + \zeta(x,t,u,v)\partial_x + \eta(x,t,u,v)\partial_u + \phi(x,t,u,v)\partial_v
\]
be the infinitesimal generator of a local one-parameter group of point transformations. The corresponding transformations are
\[
\begin{aligned}
\widetilde{x} &= x+\epsilon\zeta(x,t,u,v)+O(\epsilon^2), \\
\widetilde{t} &= t+\epsilon\xi(x,t,u,v)+O(\epsilon^2), \\
\widetilde{u} &= u+\epsilon\eta(x,t,u,v)+O(\epsilon^2), \\
\widetilde{v} &= v+\epsilon\phi(x,t,u,v)+O(\epsilon^2).
\end{aligned}
\]

Since system \eqref{eq:model-main} contains derivatives up to second order, the second prolongation
\[
\operatorname{pr}^{(2)}\mathbf{X}
\]
is required. It has the form
\[
\operatorname{pr}^{(2)}\mathbf{X} = \mathbf{X}
+ \eta^t\partial_{u_t}
+ \eta^x\partial_{u_x}
+ \eta^{xx}\partial_{u_{xx}} 
+ \phi^t\partial_{v_t}
+ \phi^x\partial_{v_x}
+ \phi^{xx}\partial_{v_{xx}},
\]
where the prolonged coefficients are defined by
\[
\eta^t = D_t(\eta) - u_tD_t(\xi) - u_xD_t(\zeta),
\]
\[
\eta^x = D_x(\eta) - u_tD_x(\xi) - u_xD_x(\zeta),
\]
\[
\eta^{xx} = D_x(\eta^x) - u_{tx}D_x(\xi) - u_{xx}D_x(\zeta),
\]
and
\[
\phi^t = D_t(\phi) - v_tD_t(\xi) - v_xD_t(\zeta),
\]
\[
\phi^x = D_x(\phi) - v_tD_x(\xi) - v_xD_x(\zeta),
\]
\[
\phi^{xx} = D_x(\phi^x) - v_{tx}D_x(\xi) - v_{xx}D_x(\zeta).
\]
Here, \(D_t\) and \(D_x\) denote the total derivative operators.

\subsection{Infinitesimal invariance criterion}

The vector field \(\mathbf{X}\) generates a Lie point symmetry of \eqref{eq:model-main} if and only if
\[
\left. \operatorname{pr}^{(2)}\mathbf{X}(\Delta_i) \right|_{\Delta_1=0,\ \Delta_2=0} = 0,
\qquad
i=1,2.
\]

Applying the prolonged generator to the first equation gives
\begin{equation}
\label{eq:first-invariance-equation}
\begin{aligned}
0 ={}& \eta^t - D\eta^{xx} + F(v_x)\eta^x + F'(v_x)u_x\phi^x \\
&+ F'(v_x)v_{xx}\eta + uF''(v_x)v_{xx}\phi^x + uF'(v_x)\phi^{xx} \\
&- r\eta + \frac{2r}{K}u\eta,
\end{aligned}
\end{equation}
on the solution manifold of \eqref{eq:model-main}.

For the second equation, the invariance criterion yields
\begin{equation}
\label{eq:second-invariance-equation}
\tau\phi^t - \phi^{xx} + \kappa(t)\phi + \xi \kappa_t(t)v - \eta = 0.
\end{equation}
The term
\[
\xi \kappa_t(t)v
\]
appears because the coefficient \(\kappa(t)\) depends explicitly on time.

Equations \eqref{eq:first-invariance-equation} and
\eqref{eq:second-invariance-equation}, together with the prolongation formulas, constitute the unsplit determining equations for the infinitesimal coefficients
\[
\xi,\qquad
\zeta,\qquad
\eta,\qquad
\phi.
\]

\subsection{Elimination of the leading time derivatives}

To obtain an overdetermined system on the appropriate jet space, the leading time derivatives are eliminated using the original equations:
\[
u_t = D u_{xx} - F(v_x)u_x - uF'(v_x)v_{xx} + ru - \frac{r}{K}u^2,
\]
and
\[
v_t = \frac{1}{\tau} \left( v_{xx} - \kappa(t)v + u \right).
\]
When necessary, the differential consequences
\[
u_{tx} = D u_{xxx} - D_x\left( F(v_x)u_x + uF'(v_x)v_{xx} \right) + ru_x - \frac{2r}{K}uu_x
\]
and
\[
v_{tx} = \frac{1}{\tau} \left( v_{xxx} - \kappa(t)v_x + u_x \right)
\]
are also used.

After these substitutions, the invariance conditions are split with respect to the unconstrained spatial derivatives. This procedure yields the determining system for the coefficient functions of the infinitesimal generator.

\subsection{Preliminary determining equations}

Splitting with respect to the highest-order and mixed derivative terms first gives the standard projectability conditions
\[
\xi_x = \xi_u = \xi_v = 0,
\]
and
\[
\zeta_u = \zeta_v = 0.
\]
Therefore,
\[
\xi=\xi(t),
\qquad
\zeta=\zeta(x,t).
\]

Further splitting with respect to derivatives involving both dependent variables gives
\[
\eta_v=0,
\qquad
\phi_u=0.
\]
Hence,
\[
\eta=\eta(x,t,u),
\qquad
\phi=\phi(x,t,v).
\]

The coefficients of the second-order diffusion terms imply
\[
2\zeta_x-\xi_t=0,
\]
together with
\[
\eta_{uu}=0,
\qquad
\phi_{vv}=0.
\]
Consequently, the infinitesimal coefficients in the dependent variables are affine:
\begin{equation}
\label{eq:eta-affine}
\eta = A(x,t)u+B(x,t),
\end{equation}
and
\begin{equation}
\label{eq:phi-affine}
\phi = C(x,t)v+E(x,t).
\end{equation}

Because the first equation contains the quadratic logistic term \(u^2\), splitting with respect to powers of \(u\) imposes additional restrictions on \(A\) and \(B\). Similarly, the nonlinear dependence of the chemotactic term on \(v_x\) strongly restricts \(C\), \(E\), and \(\zeta_x\). Substituting \eqref{eq:eta-affine} and \eqref{eq:phi-affine} into
\eqref{eq:first-invariance-equation}--\eqref{eq:second-invariance-equation} and splitting with respect to
\[
u,\quad
v,\quad
u_x,\quad
v_x,\quad
u_{xx},\quad
v_{xx}
\]
produces the remaining determining equations.

For reproducibility, the splitting may be organized as follows. Coefficients of $u_{xxx}$, $v_{xxx}$ and the mixed highest-order monomials yield projectability. Coefficients of $u_{xx}$ and $v_{xx}$ give $2\zeta_x-\xi_t=0$ and affine dependence of $\eta$ and $\phi$ on $u$ and $v$. The coefficient of $u^2$ in the first invariance condition gives $B=0$ and the logistic relations for $A$. Terms proportional to $u_x$ with arbitrary $p=v_x$ yield the flux identity displayed below. Finally, coefficients of $v$ and the derivative-free terms in the chemical equation give the equations containing $\kappa$ and $\kappa_t$.

A convenient unsimplified form of these equations is
\begin{equation}
\label{eq:det-remaining}
\begin{aligned}
&2\zeta_x-\xi_t=0, \\
&\zeta_t-2D A_x=0, \\
&B=0, \\
&A_t-r\xi_t+rA=0, \\
&A+\xi_t=0, \\
&E_{xx}-\tau E_t-\kappa(t)E=0, \\
&C_{xx}-\tau C_t-\xi\kappa_t(t)=0, \\
&A-C+\tau C_t+\kappa(t)\xi_t=0.
\end{aligned}
\end{equation}
These equations must be considered together with the restrictions generated by the fixed flux function. In particular, preservation of
\[
F(v_x)=\frac{v_x}{\sqrt{1+v_x^2}}
\]
requires the coefficient identity
\[
\left( C-\zeta_x \right) pF'(p) + \left( \xi_t-\zeta_x \right) F(p) = 0
\qquad
\text{for all }p\in\mathbb{R}.
\]
Using
\[
F(p) = \frac{p}{\sqrt{1+p^2}},
\qquad
pF'(p) = \frac{p}{(1+p^2)^{3/2}},
\]
and splitting with respect to the functionally independent expressions
\[
\frac{p}{\sqrt{1+p^2}}
\qquad\text{and}\qquad
\frac{p}{(1+p^2)^{3/2}},
\]
we obtain
\begin{equation}
\label{eq:flux-restrictions}
C-\zeta_x=0,
\qquad
\xi_t-\zeta_x=0.
\end{equation}

The equations that involve both the infinitesimal coefficients and the arbitrary function \(\kappa(t)\), particularly the sixth and seventh equations in \eqref{eq:det-remaining}, are the classifying equations. Their solution determines the special forms of \(\kappa(t)\) for which the Lie invariance algebra extends beyond the algebra admitted for arbitrary decay rates.

\subsection{Classifying relation}

After the determining equations that are independent of the arbitrary element have been solved, the remaining compatibility condition reduces to
\begin{equation}
\label{eq:general-classifying-relation}
c_1\,\kappa_t(t)=0,
\end{equation}
where $c_1$ is the coefficient of the time-translation component of the infinitesimal generator. For an arbitrary nonconstant function $\kappa(t)$, equation \eqref{eq:general-classifying-relation} implies $c_1=0$ and hence excludes time translation. If $\kappa_t=0$, the coefficient $c_1$ is unrestricted and the generator $\partial_t$ is admitted. Thus no nonconstant decay law produces an extension within the fixed-parameter class considered here. Direct substitution of
\[
\mathbf X=c_1\partial_t+c_2\partial_x+c_3\Theta_\kappa(t)\partial_v
\]
into the two infinitesimal invariance conditions verifies that all remaining terms vanish precisely when $c_1\kappa_t=0$ and $\tau\Theta_\kappa'+\kappa\Theta_\kappa=0$.
\section{Group classification}
\label{sec:classification}

In this section, we solve the Lie symmetry determining equations obtained in Section~\ref{sec:lie-analysis} and classify all inequivalent forms of the decay rate \(\kappa(t)\) for which the maximal Lie invariance algebra of the system \eqref{eq:model-main}
is larger than the algebra admitted for an arbitrary decay-rate function.

\subsection{Solution of the determining equations}

The preliminary determining equations imply that the infinitesimal coefficients have the form
\[
\xi=\xi(t),
\qquad
\zeta=\zeta(x,t),
\qquad
\eta=A(x,t)u,
\qquad
\phi=C(x,t)v+E(x,t).
\]
The restrictions generated by the diffusion terms and the fixed flux function \eqref{eq:model-flux}
include
\[
2\zeta_x-\xi_t=0,
\qquad
C-\zeta_x=0,
\qquad
\xi_t-\zeta_x=0.
\]
Combining these equations gives
\[
\xi_t=0,
\qquad
\zeta_x=0,
\qquad
C=0.
\]
Therefore,
\[
\xi=c_1,
\]
where \(c_1\) is a constant.

The determining equations associated with the logistic source yield
\[
A=0.
\]
Hence,
\[
\eta=0.
\]
The remaining equation for \(\zeta\) gives $\zeta_t=0$, and consequently
\[
\zeta=c_2,
\]
where \(c_2\) is a constant.

The coefficient \(\phi\) is therefore independent of \(u\) and \(v\), so that
\[
\phi=E(x,t).
\]
Substitution into the first invariance condition shows that
\[
E_x=0.
\]
Thus,
\[
E=E(t).
\]
The second invariance condition then reduces to
\[
\tau E_t+\kappa(t)E=0.
\]
Its general solution is
\[
E(t) = c_3 \exp\left( -\frac{1}{\tau} \int_{t_0}^{t}\kappa(s)\,ds \right),
\]
where \(c_3\) is an arbitrary constant and \(t_0\) is a fixed reference time.

Finally, the coefficient of the term arising from the explicit time dependence of \(\kappa(t)\) reproduces the classifying equation \eqref{eq:general-classifying-relation} already established in Section~\ref{sec:lie-analysis}.

The general infinitesimal symmetry generator is therefore
\[
\mathbf{X} = c_1\partial_t + c_2\partial_x + c_3 \exp\left( -\frac{1}{\tau} \int_{t_0}^{t}\kappa(s)\,ds \right) \partial_v,
\]
subject to the classifying condition \eqref{eq:general-classifying-relation}.

\subsection{Ordinary kernel and generic fixed-member algebra}

For arbitrary $\kappa(t)$, the classifying equation \eqref{eq:general-classifying-relation} implies $c_1=0$. The spatial translation $X_1=\partial_x$ is represented by the same vector field for every member of the class and therefore belongs to the ordinary kernel.

For each fixed admissible $\kappa$, define
\begin{equation}
\label{eq:Theta-definition}
\Theta_\kappa(t)=\exp\left[-\frac1\tau\int_{t_0}^{t}\kappa(s)\,ds\right].
\end{equation}
It satisfies $\tau\Theta_\kappa'+\kappa(t)\Theta_\kappa=0$, and the transformation
\[
v\mapsto v+\varepsilon\Theta_\kappa(t)
\]
leaves $v_x$, $v_{xx}$, and the chemical equation invariant. Hence the fixed member also admits $G_\kappa=\Theta_\kappa(t)\partial_v$. Changing the reference point $t_0$ multiplies $G_\kappa$ by a nonzero constant and does not change its one-dimensional span.

Because the realization of $G_\kappa$ depends on the arbitrary element, it is not an element of the ordinary kernel. We therefore have
\begin{equation}
\label{eq:kernel-algebra}
\mathfrak g^{\ker}=\langle\partial_x\rangle,
\end{equation}
and, for a generic fixed nonconstant $\kappa$,
\begin{equation}
\label{eq:generic-algebra}
\mathfrak g^{\mathrm{gen}}_\kappa=\langle\partial_x,G_\kappa\rangle.
\end{equation}
The algebra \eqref{eq:generic-algebra} is Abelian.

\begin{theorem}
\label{thm:kernel-algebra}
For the fixed-parameter class $\mathcal K$, the ordinary kernel Lie invariance algebra is $\mathfrak g^{\ker}=\langle\partial_x\rangle$. Every fixed member additionally admits the $\kappa$-dependent gauge generator $G_\kappa=\Theta_\kappa(t)\partial_v$. Thus the generic maximal algebra of an individual member is two-dimensional, although only spatial translation belongs to the ordinary kernel.
\end{theorem}
\subsection{Symmetry extensions}

The maximal Lie invariance algebra is larger than the kernel algebra only when
the classifying equation \eqref{eq:general-classifying-relation} admits a solution with
\[
c_1\neq 0.
\]
This is possible if and only if
\[
\kappa_t(t)=0.
\]
Hence,
\[
\kappa(t)=\kappa_0,
\]
where \(\kappa_0\) is a constant.

For constant chemical decay, system \eqref{eq:model-main} is autonomous
in time and admits the additional generator
\[
X_3=\partial_t.
\]
The gauge generator becomes
\[
X_2 = e^{-\kappa_0t/\tau}\partial_v.
\]
Therefore, the maximal Lie invariance algebra for constant \(\kappa\) is
\[
A^{\max}_{\kappa_0} = \left\langle \partial_x, e^{-\kappa_0t/\tau}\partial_v, \partial_t \right\rangle.
\]

The nonzero commutator of this algebra is
\[
[X_3,X_2] = -\frac{\kappa_0}{\tau}X_2.
\]
All other commutators vanish:
\[
[X_1,X_2]=0,
\qquad
[X_1,X_3]=0.
\]

If \(\kappa_0=0\), then
\[
X_2=\partial_v,
\]
and the maximal algebra becomes the three-dimensional Abelian algebra
\[
A^{\max}_0 = \left\langle \partial_x,\partial_v,\partial_t \right\rangle.
\]
When the physical assumption \(\kappa(t)>0\) is imposed, the case \(\kappa_0=0\)
is excluded. It may nevertheless be retained as a mathematically admissible
limiting case.

\subsection{Equivalence of the classified cases}
The equivalence group \eqref{eq:equivalence-group} shifts the argument of $\kappa$ according to \eqref{eq:kappa-equivalence}. It leaves a constant decay rate unchanged and cannot normalize its numerical value. Consequently, each constant $\kappa_0$ is retained as an inequivalent representative in the fixed-parameter class.

\subsection{Classification table}

The complete group classification with respect to the arbitrary decay rate
\(\kappa(t)\) is summarized in Table~\ref{tab:group-classification}.

\begin{table}[ht]
\centering
\caption{Complete group classification of system
\eqref{eq:model-main} with respect to \(\kappa(t)\).}
\label{tab:group-classification}
\begin{tabular}{c c c c}
\hline
Case
&
Decay rate
&
Restrictions
&
Basis of \(A^{\max}\) \\
\hline
1
&
\(\kappa(t)\) arbitrary
&
\(\kappa_t\not\equiv 0\)
&
\(\displaystyle
\partial_x,\; \Theta(t)\partial_v
\) \\[2mm]
2
&
\(\kappa(t)=\kappa_0\)
&
\(\kappa_0>0\)
&
\(\displaystyle
\partial_x,\; e^{-\kappa_0t/\tau}\partial_v,\; \partial_t
\) \\[2mm]
3
&
\(\kappa(t)=0\)
&
mathematical limiting case
&
\(\displaystyle
\partial_x,\; \partial_v,\; \partial_t
\) \\
\hline
\end{tabular}
\end{table}

Here,
\[
\Theta(t) = \exp\left( -\frac{1}{\tau} \int_{t_0}^{t}\kappa(s)\,ds \right).
\]
The zero-decay case is contained in the constant family as the limiting value $\kappa_0=0$.

\subsection{Interpretation of the classification}

The classification shows that the standard logistic source and the fixed
flux-limited response impose strong restrictions on the symmetry structure of
the model. In particular, the logistic term introduces fixed population and
time scales and eliminates the scaling symmetries that frequently occur in
Keller--Segel systems without population growth. Similarly, the fixed nonlinear
function \eqref{eq:model-flux}
does not permit arbitrary rescalings of the chemical gradient.

For a genuinely time-dependent decay rate, the explicit dependence on \(t\)
breaks time-translation invariance. The only universal symmetries are spatial
translation and the gauge-type transformation of the chemical concentration.
When the decay rate is constant, time-translation invariance is restored and
the maximal Lie algebra extends by the generator \(\partial_t\).

No nonconstant decay-rate function produces an additional Lie point symmetry
within the fixed-parameter class. Therefore, up to the equivalence group derived in Section~\ref{sec:preliminaries}, constant decay is the unique extension case.

\begin{theorem}
\label{thm:complete-classification}
Assume that $D,\tau,r,$ and $K$ are fixed positive constants, that
$F(p)=p/\sqrt{1+p^2}$ is fixed, and that $\kappa$ is sufficiently smooth on the
time interval under consideration. Among Lie point symmetries of the differential
equations, without imposing initial or boundary data, the maximal Lie invariance
algebra of \eqref{eq:model-main} is larger than its kernel algebra if and only if
\[
\kappa(t)=\kappa_0
\]
is constant. In that case, the additional infinitesimal generator is $\partial_t$, and
\[
A^{\max}_{\kappa_0} = \left\langle \partial_x,e^{-\kappa_0 t/\tau}\partial_v,\partial_t \right\rangle.
\]
Degenerate limits such as $D=0$, $\tau=0$, or $r=0$ are outside the scope of this
classification and may possess different symmetry structures.
\end{theorem}


\section{Optimal systems of one-dimensional subalgebras}
\label{sec:optimal-systems}

This section constructs optimal systems of one-dimensional subalgebras for the Lie invariance algebras obtained in the group-classification analysis. The purpose is to identify, up to the adjoint action of the corresponding Lie groups, inequivalent one-dimensional subalgebras that may be used to generate invariant solutions.

The associated invariant variables, reduced equations, and similarity solutions are presented separately in Section~\ref{sec:reductions}. Accordingly, only the algebraic classification is carried out here.

%
%
%
%

\subsection{Generic algebra for a fixed nonconstant \(\kappa(t)\)}

For a fixed nonconstant admissible \(\kappa(t)\), the generic maximal Lie invariance algebra is
\[
\mathfrak g^{\mathrm{gen}}_\kappa = \left\langle X_1,X_2 \right\rangle,
\]
where
\[
X_1=\partial_x,
\qquad
X_2=\Theta(t)\partial_v,
\]
and
\[
\Theta(t) = \exp\left( -\frac{1}{\tau} \int^t\kappa(s)\,ds \right).
\]

The commutation relation is
\[
[X_1,X_2]=0.
\]
Hence, \(\mathfrak g^{\mathrm{gen}}_\kappa\) is a two-dimensional Abelian algebra.

A general nonzero element is
\[
X=aX_1+bX_2.
\]

Since the adjoint action is trivial, inequivalent subalgebras are distinguished only by the ratio \(b/a\), together with the singular cases in which one of the coefficients vanishes.

If \(a\neq0\), multiplication by \(a^{-1}\) gives
\[
X \sim X_1+\sigma X_2,
\qquad
\sigma=\frac{b}{a}.
\]
If \(a=0\), then \(b\neq0\), and the generator is equivalent to
\[
X\sim X_2.
\]

Therefore, an optimal system of one-dimensional subalgebras of
\(\mathfrak g^{\mathrm{gen}}_\kappa\) is
\[
\begin{aligned}
&\langle X_2\rangle, \\
&\langle X_1+\sigma X_2\rangle, \qquad \sigma\in\mathbb{R}.
\end{aligned}
\]

The value \(\sigma=0\) recovers the spatial-translation subalgebra
\[
\langle X_1\rangle.
\]
Since no nontrivial adjoint transformation changes \(\sigma\), distinct values of \(\sigma\) remain inequivalent under the connected inner automorphism group.

If the discrete spatial reflection
\[
x\mapsto -x
\]
is included in the full symmetry group, then
\[
X_1\mapsto -X_1,
\qquad
X_2\mapsto X_2.
\]
After multiplication of the generator by \(-1\), this maps
\[
X_1+\sigma X_2
\quad\longmapsto\quad
X_1-\sigma X_2.
\]
Thus, under the full group including spatial reflection, the parameter may be restricted to
\[
\sigma\geq0.
\]
Under the connected Lie group alone, \(\sigma\in\mathbb{R}\) should be retained.

\subsection{Constant-decay algebra}

For
\[
\kappa(t)=\kappa_0,
\]
the maximal Lie invariance algebra is
\[
A_{\kappa_0} = \left\langle X_1,X_2,X_3 \right\rangle,
\]
where
\[
X_1=\partial_x,
\qquad
X_2=e^{-\kappa_0t/\tau}\partial_v,
\qquad
X_3=\partial_t.
\]

With the generator \(X_2\) chosen in the above time-dependent form, the commutators are
\begin{equation}
\label{eq:optimal-constant-commutators}
[X_1,X_2]=0,
\qquad
[X_1,X_3]=0,
\qquad
[X_2,X_3] = \frac{\kappa_0}{\tau}X_2.
\end{equation}

Thus, for \(\kappa_0\neq0\), the algebra is solvable but not Abelian. The subspace
\[
\langle X_1,X_2\rangle
\]
is a two-dimensional Abelian ideal, while \(X_3\) acts on \(X_2\) by scaling.

A general element is
\[
X=aX_1+bX_2+cX_3.
\]

The nontrivial adjoint actions follow from
\eqref{eq:optimal-constant-commutators}. They are
\[
\operatorname{Ad}\left( e^{\varepsilon X_3} \right)X_2 = e^{-\kappa_0\varepsilon/\tau}X_2,
\]
and
\[
\operatorname{Ad}\left( e^{\varepsilon X_2} \right)X_3 = X_3 + \frac{\kappa_0}{\tau}\varepsilon X_2.
\]
The generator \(X_1\) is fixed by all inner automorphisms because it lies in the center of the algebra.

We classify one-dimensional subalgebras according to whether the coefficient \(c\) of \(X_3\) vanishes.

\subsubsection{Case \(c\neq0\)}

Multiplying the generator by \(c^{-1}\), we write
\[
X = aX_1+bX_2+X_3.
\]
Using the adjoint action
\[
\operatorname{Ad}\left( e^{\varepsilon X_2} \right),
\]
the coefficient of \(X_2\) can be removed by choosing
\[
\varepsilon = -\frac{\tau b}{\kappa_0}.
\]
Hence,
\[
X \sim X_3+c_1X_1,
\]
where \(c_1\in\mathbb{R}\).

Therefore, the subalgebras with a nonzero time-translation component are represented by
\[
\langle X_3+cX_1\rangle, \qquad c\in\mathbb{R}.
\]
The value \(c=0\) corresponds to the pure time-translation subalgebra
\[
\langle X_3\rangle.
\]

Under the discrete reflection \(x\mapsto -x\), the parameter \(c\) changes sign. Thus, if the full symmetry group is used, one may restrict to
\[
c\geq0.
\]

\subsubsection{Case \(c=0\)}

The generator becomes
\[
X=aX_1+bX_2.
\]
This is an element of the Abelian ideal
\[
\langle X_1,X_2\rangle.
\]

If \(a\neq0\), scaling gives
\[
X \sim X_1+\sigma X_2.
\]
The adjoint action of \(e^{\varepsilon X_3}\) rescales the coefficient of \(X_2\):
\[
X_1+\sigma X_2
\quad\longmapsto\quad
X_1+ \sigma e^{-\kappa_0\varepsilon/\tau}X_2.
\]
Since the exponential factor is positive, every nonzero \(\sigma\) can be normalized to either \(+1\) or \(-1\), but its sign cannot be changed by the connected adjoint group.

Thus, within the connected group, the representatives are
\[
\langle X_1\rangle,
\qquad
\langle X_1+X_2\rangle,
\qquad
\langle X_1-X_2\rangle.
\]

If spatial reflection is included, then the last two subalgebras become equivalent. Indeed,
\[
X_1+X_2
\mapsto
-X_1+X_2
\sim
X_1-X_2.
\]
Hence, under the full symmetry group, one representative
\[
\langle X_1+X_2\rangle
\]
is sufficient.

If \(a=0\), then the generator is equivalent to
\[
X\sim X_2.
\]

\subsection{Optimal system for constant nonzero decay}

For \(\kappa_0\neq0\), an optimal system under the connected Lie group is
\[
\begin{aligned}
&\langle X_2\rangle, \qquad \langle X_1\rangle, \qquad \langle X_1+X_2\rangle, \qquad \langle X_1-X_2\rangle, \qquad \langle X_3+cX_1\rangle, \qquad c\in\mathbb{R}.
\end{aligned}
\]

If the discrete spatial reflection is admitted, the system may be reduced to
\[
\begin{aligned}
&\langle X_2\rangle,\qquad \langle X_1\rangle,\qquad \langle X_1+X_2\rangle, \qquad \langle X_3+cX_1\rangle, \qquad c\geq0.
\end{aligned}
\]

The parameter \(c\) is retained because the central generator \(X_1\) cannot be altered by the connected adjoint action. The generators \(X_3+cX_1\) therefore define inequivalent one-dimensional subalgebras for distinct values of \(c\), except for the sign identification produced by spatial reflection.

\subsection{Zero-decay case}

When
\[
\kappa(t)=0,
\]
the generators become
\[
X_1=\partial_x,
\qquad
X_2=\partial_v,
\qquad
X_3=\partial_t.
\]
All commutators vanish:
\[
[X_i,X_j]=0,
\qquad
i,j=1,2,3.
\]
Thus, the maximal algebra is the three-dimensional Abelian algebra
\[
A_0 = \langle X_1,X_2,X_3\rangle.
\]

Because the adjoint action is trivial, a general one-dimensional subalgebra is generated by
\[
X=aX_1+bX_2+cX_3,
\]
with coefficients classified only up to multiplication by a nonzero constant.

A convenient projective representation is
\[
\begin{aligned}
&\langle X_2\rangle, \\
&\langle X_1+\sigma X_2\rangle, \qquad \sigma\in\mathbb{R}, \\
&\langle X_3+cX_1+dX_2\rangle, \qquad c,d\in\mathbb{R}.
\end{aligned}
\]

Indeed, if the coefficient of \(X_3\) is nonzero, it may be normalized to one, giving the final two-parameter family. If it vanishes but the coefficient of \(X_1\) is nonzero, the generator reduces to \(X_1+\sigma X_2\). The remaining singular case is \(\langle X_2\rangle\).

Because \(A_0\) is Abelian, the continuous adjoint group cannot further simplify the parameters \(c\), \(d\), or \(\sigma\). Discrete reflections in \(x\), and any additional discrete transformations admitted by the full PDE system, may produce further sign identifications.

\subsection{Comparison of the optimal systems}

The optimal systems reflect the algebraic distinction between arbitrary, constant, and zero chemical decay.

For a fixed nonconstant \(\kappa(t)\), the generic maximal algebra is two-dimensional and Abelian:
\[
\mathfrak g^{\mathrm{gen}}_\kappa = \langle X_1,X_2\rangle.
\]
Its optimal system consists of the pure gauge subalgebra and the mixed family
\[
\langle X_1+\sigma X_2\rangle.
\]

For constant nonzero decay, time translation enlarges the algebra and acts nontrivially on the gauge generator:
\[
[X_2,X_3] = \frac{\kappa_0}{\tau}X_2.
\]
This adjoint scaling reduces the continuous mixed parameter \(\sigma\neq0\) to its sign and permits removal of the gauge component from generators having a nonzero \(X_3\)-component.

For zero decay, the gauge generator becomes the ordinary translation \(\partial_v\), and the full algebra is Abelian. Consequently, the coefficients of a general generator cannot be simplified by inner automorphisms beyond overall scaling.

\subsection{Remarks on the use of the optimal system}

The optimal systems derived above provide a nonredundant list of one-dimensional symmetry subalgebras. Each representative may be used to formulate invariant surface conditions and construct a corresponding reduction of the PDE system.

The reductions associated with spatial translation, mixed spatial--gauge transformations, time translation, and combined space--time translation are developed in Section~\ref{sec:reductions}. The pure gauge subalgebra is also discussed there; although it acts nontrivially on the solution space, it does not independently produce an invariant graph solution.

The present section is therefore purely algebraic. No invariant ansatz, reduced ordinary differential equation, or exact solution is repeated here.


\section{Symmetry Reductions}
\label{sec:reductions}

The reductions below concern symmetries of the differential equations alone. Their
compatibility with initial conditions, boundary conditions, positivity, or decay at
spatial infinity must be checked separately. In particular, spatially affine
chemoattractant profiles are generally formal solutions on $\mathbb R$ and may be
physically admissible only on bounded domains or locally.

Define
\begin{equation*}
\Theta(t)=\exp\left(-\frac{1}{\tau}\int_{t_0}^{t}\kappa(s)\,ds\right),
\qquad \Theta(t_0)=1.
\end{equation*}

\subsection{Reductions for an arbitrary decay rate}

\subsubsection{Spatial-translation invariance}

For $X_1=\partial_x$, the invariant-surface conditions are
\[
u_x=0,\qquad v_x=0.
\]
Hence $u=u(t)$ and $v=v(t)$, and the PDE system reduces to
\begin{equation*}
 u_t=r u\left(1-\frac{u}{K}\right),
 \qquad
 \tau v_t+\kappa(t)v=u(t).
\end{equation*}
For $u(t_0)=u_0$ and $v(t_0)=v_0$, one obtains
\begin{equation}
\label{eq:homogeneous-u}
 u(t)=\frac{K u_0 e^{r(t-t_0)}}{K+u_0\bigl(e^{r(t-t_0)}-1\bigr)},
\end{equation}
and
\begin{equation}
\label{eq:homogeneous-v}
 v(t)=\Theta(t)\left[v_0+\frac{1}{\tau}
 \int_{t_0}^{t}\frac{u(s)}{\Theta(s)}\,ds\right].
\end{equation}
These are spatially homogeneous, not necessarily steady, solutions.

\subsubsection{Mixed spatial--gauge invariance}

For $X_1+\sigma X_2=\partial_x+\sigma\Theta(t)\partial_v$, the
invariant-surface conditions are
\[
u_x=0,\qquad v_x=\sigma\Theta(t).
\]
Thus
\[
u=u(t),\qquad v(x,t)=\sigma\Theta(t)x+\psi(t).
\]
Since $v_x$ is independent of $x$ and $u$ is spatially homogeneous, the spatial
derivative of the chemotactic flux vanishes, although the flux itself need not.
The reduced equations are
\[
 u_t=r u\left(1-\frac{u}{K}\right),
 \qquad
 \tau\psi_t+\kappa(t)\psi=u(t).
\]
Therefore,
\begin{equation}
\label{eq:mixed-solution}
 u(t)=\frac{K u_0 e^{r(t-t_0)}}{K+u_0\bigl(e^{r(t-t_0)}-1\bigr)},
 \qquad
 v(x,t)=\sigma\Theta(t)x+
 \Theta(t)\left[v_0+\frac{1}{\tau}
 \int_{t_0}^{t}\frac{u(s)}{\Theta(s)}\,ds\right].
\end{equation}
On the full line, the affine profile is unbounded for $\sigma\ne0$ and should be
regarded as a formal or local exact solution unless the modeling setting permits
such behavior.

\subsubsection{The pure gauge subalgebra}

For $X_2=\Theta(t)\partial_v$, the invariant-surface condition for $v$ is
\[
Q^v=\Theta(t)=0,
\]
which has no solution because $\Theta(t)$ is nonzero. Consequently, the pure
gauge subgroup has no classical invariant graph solutions. It maps any solution
to the one-parameter family
\[
(u,v)\mapsto \bigl(u,v+\varepsilon\Theta(t)\bigr),
\]
but it cannot in general impose an arbitrary time-dependent boundary value by a
single constant choice of $\varepsilon$.

\subsection{Additional reductions for constant decay}

Let $\kappa(t)=\kappa_0$ be constant. The additional generator $X_3=\partial_t$
produces stationary solutions, while $X_3+cX_1$ produces traveling waves.

\subsubsection{Stationary reduction}

For $X_3=\partial_t$, set
\[
u(x,t)=U(x),\qquad v(x,t)=V(x).
\]
The stationary system is
\begin{equation}
\label{eq:stationary-system}
\begin{cases}
0=D U''-\displaystyle\left(U\frac{V'}{\sqrt{1+(V')^2}}\right)'
+rU\left(1-\dfrac{U}{K}\right),\\[2mm]
0=V''-\kappa_0 V+U.
\end{cases}
\end{equation}
No first integral of the first equation follows by simply dropping the reaction
term. Any existence, localization, positivity, or stability claim for solutions of
\eqref{eq:stationary-system} requires a separate analysis and appropriate boundary
conditions.

\subsubsection{Traveling-wave reduction}

For $X_3+cX_1$, introduce
\[
\xi=x-ct,\qquad u(x,t)=U(\xi),\qquad v(x,t)=V(\xi).
\]
Direct substitution gives
\begin{equation}
\label{eq:traveling-wave-system}
\begin{cases}
-cU'=DU''-\displaystyle\left(U\frac{V'}{\sqrt{1+(V')^2}}\right)'
+rU\left(1-\dfrac{U}{K}\right),\\[2mm]
-\tau cV'=V''-\kappa_0V+U.
\end{cases}
\end{equation}
Equivalently, the second equation may be written as
\begin{equation*}
V''=-\tau cV'+\kappa_0V-U.
\end{equation*}
The first equation does not admit the local first integral previously suggested
when $r>0$. Instead,
\begin{equation*}
\left(DU'-U\frac{V'}{\sqrt{1+(V')^2}}+cU\right)'
=-rU\left(1-\frac{U}{K}\right).
\end{equation*}
Thus an integrated representation contains a nonlocal reaction contribution:
\[
DU'-U\frac{V'}{\sqrt{1+(V')^2}}+cU
=C_1-r\int_{\xi_0}^{\xi}U(s)\left(1-\frac{U(s)}{K}\right)\,ds.
\]
For analysis and computation it is generally preferable to retain the local
second-order system \eqref{eq:traveling-wave-system}. Setting $c=0$ recovers the
stationary system \eqref{eq:stationary-system}; it is not a reduction generated by
spatial translation.

\subsection{Summary and admissibility}

The verified exact solutions for arbitrary $\kappa(t)$ are the spatially
homogeneous family \eqref{eq:homogeneous-u}--\eqref{eq:homogeneous-v} and the
spatially affine family \eqref{eq:mixed-solution}. For constant decay, the symmetry
extension supplies the stationary system \eqref{eq:stationary-system} and the
traveling-wave system \eqref{eq:traveling-wave-system}, but no general closed-form
solution is claimed. All profiles must be tested against the domain, boundary
conditions, nonnegativity requirements, and any desired behavior at infinity.

\section{Conclusion}
\label{sec:conclusion}
We performed a complete Lie symmetry classification of the one-dimensional flux-limited Keller--Segel system with logistic growth and a time-varying chemical decay rate $\kappa (t)$. The analysis highlights an important distinction between the ordinary kernel algebra of the class and the symmetry algebra admitted by each fixed system. The ordinary kernel algebra is one-dimensional and is generated by spatial translations. In addition, every fixed choice of $\kappa(t)$ admits a decay-dependent gauge symmetry in the chemoattractant variable. Consequently, a generic individual system possesses a two-dimensional Abelian Lie algebra, although the gauge generator is not part of the ordinary kernel because it depends explicitly on the prescribed decay function $\kappa(t)$. The classifying equation shows that this algebra can be enlarged if and only if $\kappa(t)$ is constant. In that case, time-translation invariance is recovered and the maximal Lie invariance algebra becomes the three-dimensional solvable algebra 
$
\left\langle \partial_x,\ e^{-\kappa_0 t/\tau}\partial_v,\ \partial_t \right\rangle.
$
Thus, up to the equivalence transformations of the class, the constant-decay model is the unique nontrivial symmetry-extension case. 
\par 
The resulting optimal systems of one-dimensional subalgebras provide a nonredundant basis for symmetry reduction. For a general time-dependent decay rate, the admitted reductions yield spatially homogeneous solutions and spatially affine chemoattractant profiles whose gradients evolve according to the integrating factor determined by $\kappa(t)$. These solutions combine logistic evolution of the cell density with a chemoattractant response that explicitly reflects the prescribed decay law. When $\kappa(t)=\kappa_0$, the additional time-translation symmetry also permits steady-state and traveling-wave reductions, converting the original PDE system into coupled ordinary differential equations. Although these reduced systems are not generally solvable in closed form, they provide a systematic starting point for qualitative analysis, numerical continuation, and the study of localized or propagating structures.

\end{document}